\documentclass[12pt]{article}
\usepackage[dvips]{graphicx}
\usepackage{amsmath,amssymb,amsthm}
\def\Bt{{\bf t}}
\def\Bx{{\bf x}}
\def\By{{\bf y}}
\def\Bz{{\bf z}}

\def\Bb{{\bf b}}

\def\BA{{\bf A}}
\def\BB{{\bf B}}
\def\BC{{\bf C}}
\def\BD{{\bf D}}

\def\BE{{\bf E}}

\def\BI{{\bf I}}
\def\Bi{{\bf i}}
\def\Bj{{\bf j}}

\def\Bt{{\bf t}}

\def\Bbeta{{\bf \beta}}

\def\Bone{{\bf 1}}
\def\Bzero{{\bf 0}}

\newtheorem{theorem}{Theorem}
\newtheorem{proposition}{Proposition}

\newtheorem{corollary}{Corollary}

\newtheorem{lemma}{Lemma}

\theoremstyle{definition}
\newtheorem{example}{Example}
\newtheorem{definition}{Definition}

\title{Minimal Markov basis for tests of main effect models for $2^{p-1}$ 
fractional factorial designs of resolution $p$}

\author{Satoshi Aoki%
\thanks{Graduate School of Science and Engineering (Science Course), Kagoshima University.}%
\ \thanks{JST, CREST.}
}

\begin{document}
\maketitle

\begin{abstract}
We consider conditional exact tests of factor effects in
design of experiments for discrete response variables.  Similarly to
the analysis of contingency tables, Markov chain Monte Carlo methods
can be used to perform exact tests, especially when large-sample
approximations of the null distributions are poor and 
the enumeration of the conditional sample space is infeasible.  
In order to construct a connected Markov chain over the 
appropriate sample space,  one approach is to compute a
Markov basis. Theoretically, a Markov basis can be characterized as a
generator of a well-specified toric ideal in a polynomial ring and is 
computed by computational algebraic software. However, the computation
of a Markov basis sometimes becomes infeasible, even for problems of
moderate sizes. In the present paper, we obtain the closed-form expression of minimal Markov bases for the main effect models of $2^{p-1}$ fractional factorial designs of resolution $p$.
\end{abstract}

\section{Introduction}
\label{sec:intro}
In the past decade, various new applications of computational
algebraic techniques to statistical problems have been developed. 
One of the first studies in this field, i.e.,
{\it computational algebraic statistics}, was conducted by Diaconis and Sturmfels
(\cite{Diaconis-Sturmfels-1998}), who introduced 
the notion of a {\it Markov basis} and presented a procedure for sampling
from discrete conditional distributions by constructing a connected, aperiodic,
and reversible Markov chain on a given sample space. Since then, a number of studies 
by both algebraists and statisticians have been published on the topic of Markov bases. 
Some results on the structure of Markov bases for
various statistical models are presented in \cite{Aoki-Hara-Takemura-2012}. 

The present paper is based on two studies in this field, 
\cite{Aoki-Takemura-2010} and \cite{Aoki-Takemura-2009b}, which
 consider a Markov basis in the context of design of experiments. 
In these studies, Markov chain Monte Carlo methods for testing
factor effects are discussed, when observations are discrete and are given 
in the two-level
or three-level regular fractional factorial designs. 
As one of the contributions of these studies, the relationship between 
the statistical models for the regular fractional factorial designs and 
contingency tables is considered through Markov bases. 
As a consequence, in order to investigate the Markov bases arising in 
problems associated with design of experiments, we can refer to published 
results on corresponding models for the contingency tables.  
For example, we see that the Markov basis for 
the main effect models of the regular $2^{5-2}_{{\rm III}}$ 
design given by the defining relation $\BA\BB\BD = \BA\BC\BE = \BI$ is 
constructed only by the
square-free degree-2 elements (\cite{Aoki-Takemura-2010}). 
This is because the corresponding model 
in the contingency tables is the conditional independence model in 
the $2\times 2\times 2$ table, 
i.e., the model that the two factors are 
independent in the $2\times 2$ tables for each level of the third factor. 
Note that the conditional independence model in the three-way contingency 
table is an example of decomposable models and we know that a 
minimal Markov basis for this class of models can be 
constructed only by square-free degree-2 elements. 
See \cite{Dobra-2003} for details. On the other hand, the structure 
of the Markov bases for the models that do not correspond 
to the decomposable models is not known in general.

In the present paper, following the Markov chain Monte Carlo approach in 
the design of experiments by \cite{Aoki-Takemura-2010} and 
\cite{Aoki-Takemura-2009b}, we present new results on the structure of 
the minimal Markov basis for the main effect models of $2^{p-1}$ fractional 
factorial designs of resolution $p$ for $p$ two-level factors. 
We herein demonstrate that the minimal Markov basis for this case 
can be constructed by square-free degree-2 elements.

The organisation of this paper is as follows.
In Section \ref{sec:mcmc}, we review the Markov chain Monte Carlo approach
for testing the fit of log-linear models when observations are 
discrete random variables. 
 The most important results are presented in 
Section \ref{sec:minimal-MB-2^p-1}.

\section{Markov chain Monte Carlo method for regular two-level
fractional factorial designs}
\label{sec:mcmc}
In this section, we introduce Markov chain Monte Carlo methods for 
testing the fit of log-linear models for regular two-level
fractional factorial designs with count observations.
We consider the design of two-level $p$ factors.
The observations are denoted as $\By = (y_1,\ldots,y_k)'$, 
where $k$ is the run size and $'$ denotes the transpose. 
We assume that the observations are counts of some events.
For simplicity, we also assume that only one observation is obtained 
for each run. 
This assumption is natural for setting the Poisson sampling scheme, 
because the set of totals for each run makes up the sufficient statistics 
for the parameters.
Write the $k\times p$ design matrix $D = (d_{ij})$, where
$d_{ij} \in \{-1,1\}$ is the level of the $j$-th factor in the $i$-th run
for $i=1,\ldots,k, j=1,\ldots,p$.
As an example of $2^{4-1}_{{\rm IV}}$ fractional factorial design, 
the design matrix
$D$ is given as
\begin{equation}
D=\left(\begin{array}{cccccccc}
+1 & +1 & +1 & +1 & -1 & -1 & -1 & -1\\
+1 & +1 & -1 & -1 & +1 & +1 & -1 & -1\\
+1 & -1 & +1 & -1 & +1 & -1 & +1 & -1\\
+1 & -1 & -1 & +1 & -1 & +1 & +1 & -1
\end{array}
\right)'.
\label{eqn:D:2^4-1}
\end{equation}
For count observations, 
it is natural to consider a Poisson distribution as a 
sampling model 
in the framework of generalized linear models.
 See \cite{McCullagh-Nelder-1989}.
The observations $\By$ are realizations from $k$ Poisson random
variables $Y_1,\ldots,Y_k$, which are mutually independently distributed
with mean parameter $\mu_i = E(Y_i), i = 1,\ldots,k$. 
We herein refer to the log-linear model given by 
\begin{equation}
\log \mu_i = \beta_0 + \beta_i d_{i1} + \cdots + \beta_p d_{ip},\ 
i=1,\ldots,k
\label{eqn:log-linear-null-model}
\end{equation}
as the {\it main effect model}.
 In matrix form, the main effect model is written as follows:
\[
\left(
\begin{array}{c}
\log\mu_1\\
\vdots\\
\log\mu_k
\end{array}
\right) = M\Bbeta,
\]
where $\Bbeta = (\beta_0,\beta_1,\ldots,\beta_p)'$ and 
\begin{equation}
M = \left(\begin{array}{cccc}
1 &  & & \\
\vdots &  & D & \\
1 &  & & 
\end{array}
\right).
\label{eqn:def-M-by-D}
\end{equation}
We refer to the $k\times (p+1)$ matrix $M$ 
as a {\it model matrix} of the main effect model.
The interpretation of the parameter $\beta_j$ 
in (\ref{eqn:log-linear-null-model}) 
is the parameter contrast for
the main effect of the $j$-th factor.
In \cite{Aoki-Takemura-2010}, 
models including various interaction effects are also considered.
In the present paper, we focus on the main effect model.

In order to judge the fitting of the main effect 
model (\ref{eqn:log-linear-null-model}), we can perform various 
goodness-of-fit tests, in which the main effect model
(\ref{eqn:log-linear-null-model}) is treated as the null model, whereas the 
saturated model is treated as the alternative model. 
Under the null model  (\ref{eqn:log-linear-null-model}), 
$\Bbeta$ is the nuisance parameter, and 
the  sufficient statistic for $\Bbeta$
is given by 
$M'\By = (\sum_{i=1}^k y_i, \sum_{i = 1}^k d_{i1}y_i, \ldots, \allowbreak\sum_{i = 1}^k d_{ip}y_i)'$. Then, the 
conditional distribution of $\By$ given the sufficient 
statistic is written as
\begin{equation}
f(\By\ |\ M'\By = M'\By^o) = \frac{1}{C(M'\By^o)} \displaystyle\prod_{i = 1}^k\frac{1}{y_i!},
\label{eqn:poisson-conditional-distribution}
\end{equation}
where $\By^o$ is the  observed count vector, and 
$C(M'\By^o)$ is the normalizing constant determined from
$M'\By^o$ written as
\begin{equation}
 C(M'\By^o) = \displaystyle\sum_{\By \in {\cal F}(M'\By^o)}\left(
\displaystyle\prod_{i = 1}^k\frac{1}{y_i!}
\right),
\label{eqn:poisson-constant}
\end{equation}
and
\begin{equation}
 {\cal F}(M'\By^o) = \{\By\ |\ M'\By = M'\By^o,\ y_i \ \mbox{is a
  nonnegative integer for}\ i = 1,\ldots,k\}.
\label{eqn:poisson-fiber}
\end{equation}
 We call (\ref{eqn:poisson-fiber}) as a fiber in Section 3.1.
Note that, by sufficiency, the conditional distribution 
does not depend on the values of the nuisance parameters. 

 The goodness-of-fit tests we consider are based on
the conditional distribution
(\ref{eqn:poisson-conditional-distribution}). There are several ways to
choose the test statistics. For example, the likelihood ratio statistic
\begin{equation}
\label{eq:g2}
T(\By) = 2\sum_{i = 1}^{k}y_i\log\frac{y_i}{\hat{\mu_i}}
\end{equation}
is frequently used, 
where $\hat{\mu_i}$ is the maximum likelihood estimate for $\mu_i$
under the null model (i.e., fitted value). 
Note that the traditional asymptotic test evaluates the upper probability
for the observed value $T(\By^o)$ based on the asymptotic distribution
$\chi_{k - p - 1}^2$. However, since the fitting of the asymptotic 
approximation may sometimes be poor, we consider Markov chain Monte Carlo 
methods by which to evaluate the $p$ values. 
Using the conditional distribution
(\ref{eqn:poisson-conditional-distribution}), the exact $p$ value is written as
\begin{equation}
 p = \displaystyle\sum_{\By \in {\cal F}(M'\By^o)}f(\By\ |\ M'\By =
 M'\By^o)\Bone(T(\By) \geq T(\By^o)),
\label{eqn:exact-p-value}
\end{equation}
where
\begin{equation}
 \Bone(T(\By) \geq T(\By^o)) 
= \left\{\begin{array}{ll}
1, & \mbox{if}\ T(\By) \geq T(\By^o),\\
0, & \mbox{otherwise}.
\end{array}
\right.
\label{eqn:exact-p-value-1}
\end{equation}
 Ideally, we would like to compute the exact
$p$ value of (\ref{eqn:exact-p-value})
and (\ref{eqn:exact-p-value-1}).
Unfortunately, however, an enumeration of all the elements in ${\cal
F}(M'\By^o)$, and hence the calculation of the normalizing constant 
$C(M'\By^o)$, is usually computationally infeasible for a large sample
space. In such situations, various Monte  Carlo methods are effective. 
In particular, we consider a Markov chain Monte Carlo method.
 Note that one of the important advantages of the Markov chain 
Monte Carlo method is that 
we  do not need to calculate the normalizing constant 
(\ref{eqn:poisson-constant}) in order to evaluate $p$ values.

In order to perform the Markov chain Monte Carlo procedure, 
we have to construct a connected, aperiodic, and reversible 
Markov chain over the 
conditional sample space (\ref{eqn:poisson-fiber})
with the stationary distribution (\ref{eqn:poisson-conditional-distribution}).
If such a chain is constructed, we can sample  from the chain
$\By^{(1)},\ldots,\By^{(T)}$ 
after discarding some initial burn-in steps 
and evaluate  the $p$ value as
\[
\hat{p} = \frac{1}{T}\sum_{t=1}^{T}\Bone(T(\By^{(t)}) \geq T(\By^o)). 
\] 
Such a chain can be easily constructed by a {\it Markov basis}. 
Once a Markov basis is
calculated, we can construct a connected, aperiodic, and reversible
Markov chain over the space (\ref{eqn:poisson-fiber}), which can be modified
so that the stationary distribution is the conditional distribution
(\ref{eqn:poisson-conditional-distribution})
by the Metropolis-Hastings procedure. See
\cite{Diaconis-Sturmfels-1998} and \cite{Hastings-1970} for details.
We present a formal definition of the Markov basis in the next section. 

The Markov basis is characterized algebraically as follows. 
Write indeterminates
$u_1,\ldots,u_k$ and consider polynomial ring $K[u_1,\ldots,u_k]$ for some
field $K$. Consider the integer kernel of the transpose of the model
matrix $M$, ${\rm Ker}_{\mathbb{Z}}M'$, 
 where 
we denote the set of integers by $\mathbb{Z} = \{0,\pm 1, \pm 2,\ldots,\}$.
For each $\Bb = (b_1,\ldots,b_k)' \in {\rm Ker}_{\mathbb{Z}}M'$, define
 a binomial in $K[u_1,\ldots,u_k]$ as
\[
f_{\Bb} = \prod_{b_j > 0}u_j^{b_j} - \prod_{b_j < 0}u_j^{-b_j}. 
\]
Then, the binomial ideal in $K[u_1,\ldots,u_k]$, 
\[
I(M') = \left<
\{f_{\Bb}\ |\ \Bb \in {\rm Ker}_{\mathbb{Z}}M'\}
\right>,
\]
is referred to as a toric ideal with the configuration $M'$. 
Let $\{f_{\Bb^{(1)}},\ldots,f_{\Bb^{(s)}}\}$ be any generating set of $I(M')$.
Then, the set of integer vectors $\{\Bb^{(1)},\ldots,\Bb^{(s)}\}$ constitutes
a Markov basis. See \cite{Diaconis-Sturmfels-1998} for details. 
\begin{example}
Consider the main effect model of the $2^{4-1}_{{\rm IV}}$ fractional
factorial design. The design matrix $D$ is given as (\ref{eqn:D:2^4-1})
and the model matrix is given as (\ref{eqn:def-M-by-D}). Because there are
$8$ runs in this design, we write indeterminates $u_1,\ldots,u_8$ and consider
polynomial ring $K[u_1,\ldots,u_8]$. The toric ideal with the configuration
$M'$ is calculated as
\[
I(M') = \left<
u_1u_8-u_2u_7,\ u_1u_8-u_3u_6,\ u_1u_8-u_4u_5
\right>. 
\]
Then we have a Markov basis
\[
\{(1,-1,0,0,0,0,-1,1),\ 
(1,0,-1,0,0,-1,0,1),\ 
(1,0,0,-1,-1,0,0,1)\}.
\]
\end{example}
In order to compute a Markov basis for  a given configuration $M'$, 
we can rely on
various algebraic software, such as 4ti2 (\cite{4ti2}). 
However, the computation of a Markov basis sometimes becomes infeasible,
even for problems of moderate size. In the present paper, we present the
closed-form expression of minimal Markov bases for the main effect models
of $2^{p-1}$ fractional factorial designs. 

\section{Minimal Markov basis for the main effect models 
of $2^{p-1}$ fractional factorial designs of resolution $p$}
\label{sec:minimal-MB-2^p-1}.

\subsection{Markov basis}
As we have demonstrated in the previous section, 
if we can construct a connected 
Markov chain over ${\cal F}(M'\By^o)$ for the observation $\By^o$, 
we can estimate the conditional $p$ value and judge the fitting of the 
main effect model. One common approach to constructing such a chain 
is to calculate a Markov basis for $M'$ in advance. 
Here we give  a formal definition
 of the Markov basis and  the necessarily notation. 
 For details, see \cite{Aoki-Hara-Takemura-2012}.
We denote the set of integers and nonnegative integers 
by $\mathbb{Z} = \{0,\pm 1,\pm 2,\ldots,\}$ 
and $\mathbb{N} = \{0,1,2,\ldots,\}$, respectively. 
We also denote the set of
$n$-dimensional vectors of elements from $\mathbb{Z}$ and $\mathbb{N}$ by 
$\mathbb{Z}^n$ and $\mathbb{N}^n$, respectively.

We refer to the set ${\cal F}(M'\By^o)$ of (\ref{eqn:poisson-fiber}) 
as a {\it fiber}.  Fiber ${\cal F}(M'\By^o)$ 
is the set of the nonnegative integer vectors $\By$ 
that has the same value of the 
sufficient statistics $M'\By$ to the observed vector $\By^o$.
 We refer to $\Bz \in \mathbb{Z}^k$ as a {\it move} for $M'$ 
if $M'\Bz = \Bzero$. If $M'\Bx = \Bt$ and $\Bz$ is a move for $M'$, we have
$M'(\Bx + \Bz) = M'\Bx = \Bt$. Therefore, by adding $\Bz$ to $\Bx$, we remain
in the same fiber as long as $\Bx + \Bz$ does not contain a negative
element. This is the reason that $\Bz$ is referred to as a move. 
Assume that we are given a finite set of moves ${\cal B}$. 
For each $\By \in \mathbb{N}^k$, we consider an undirected graph
$G(\By^o, {\cal B})$, the vertices of which are the elements of a fiber 
${\cal F}(M'\By^o)$. We draw an undirected edge 
between $\Bx$ and $\By$ in ${\cal F}(M'\By^o)$  
if there exists $\Bz \in {\cal B}$ such that $\By = \Bx + \Bz$ or 
$\By = \Bx - \Bz$. The Markov basis is defined as follows.
\begin{definition}
A finite set ${\cal B}$ of moves for $M'$ is
called a {\it Markov basis} for $M'$  if $G(\By^o, {\cal B})$ is 
connected for every $\By^o \in \mathbb{N}^k$.
\end{definition}
Assume that ${\cal B}$ is a Markov basis for $M'$. Then, for 
 any 
$\Bx, \By \in {\cal F}(M'\By^o)$ with $\Bx \neq \By$  
and  for any $\By^o \in \mathbb{N}^k$, 
there exists a sequence of moves $\Bz_{1}, \ldots, \Bz_{L}$ from
$\cal B$ and $\epsilon_j=\pm 1, \ j=1,\ldots,L$, such that
$\By= \Bx + \sum_{j=1}^L \epsilon_j \Bz_{j}$ and
\[
 \Bx + \sum_{j=1}^h \epsilon_j \Bz_{j} \in {\cal F}(M'\By^o), \qquad
 h=1,\ldots,L-1,
\]
i.e., we can move from $\Bx$ to $\By$ 
 using moves in 
${\cal B}$ without
generating negative elements along the way.    
In this case, we say that $\By$ is {\it accessible} from $\Bx$ by ${\cal B}$
and denote this by $\Bx \sim \By\ ({\rm mod}\ {\cal B})$. 
Obviously, the notion of
accessibility is symmetric and transitive.
Allowing moves to be $\Bzero$ also yields reflexivity. 
Therefore, accessibility by ${\cal B}$ is an equivalence relationship, and
each fiber ${\cal F}(M'\By^o)$ is partitioned into disjoint
equivalence classes by moves of ${\cal B}$. We refer to these equivalence
classes as ${\cal B}$-equivalence classes of ${\cal F}(M'\By^o)$.
Since the notion of accessibility is symmetric, we also say that
$\Bx$ and $\By$ are mutually accessible
by ${\cal B}$ if
$\Bx \sim \By \ ({\rm mod\ } {\cal B})$. 

A Markov basis ${\cal B}$ is {\it minimal} if no proper subset of ${\cal B}$
is a Markov basis. A minimal Markov basis always exists, because
from any Markov basis, we can remove redundant elements one by one
until none of the  remaining elements can be removed.
For the structure of the minimal Markov basis, see Chapter 5 
of \cite{Aoki-Hara-Takemura-2012}. 

\subsection{Main result}
Next, we present the fundamental results of the present paper. We consider the 
main effect model of $2^{p-1}$ fractional
factorial design of resolution $p$ for $p$ two-level factors.  
The $2^{p-1}$ fractional factorial design of resolution $p$ has a clear
advantage over other regular $1/2$-fractional factorial designs because 
it possesses the largest resolution. For example, $2^{5-1}_{{\rm V}}$ 
fractional factorial design is the only regular $1/2$-fractional factorial 
design where all two-factor interactions are separately estimable if we
ignore three-factor interactions.  For 
$2^{p-1}$ fractional
factorial design of resolution $p$, 
the run size is $k = 2^{p-1}$, and the model matrix $M$ is 
$k\times (p+1) = 2^{p-1} \times (p+1)$. For example, $M$ for $p = 4$ is given by
\begin{equation}
M=\left(\begin{array}{cccccccc}
+1 & +1 & +1 & +1 & +1 & +1 & +1 & +1\\
+1 & +1 & +1 & +1 & -1 & -1 & -1 & -1\\
+1 & +1 & -1 & -1 & +1 & +1 & -1 & -1\\
+1 & -1 & +1 & -1 & +1 & -1 & +1 & -1\\
+1 & -1 & -1 & +1 & -1 & +1 & +1 & -1
\end{array}
\right)'.
\label{eqn:M-p=4}
\end{equation}
A move $\Bz = (z_1,\ldots,z_k)'$ 
for $M'$ is referred to as a square-free degree-$2$ move if 
\[
z_i = \left\{\begin{array}{cl}
+1, & \mbox{for}\ i = i_1,i_2,\\ 
-1, & \mbox{for}\ i = i_3,i_4,\\
0, & \mbox{otherwise} 
\end{array}
\right.
\]
for some distinct $i_1,i_2,i_3,i_4 \in \{1,\ldots,k\}$. 
Our main theorem is as follows.
\begin{theorem}
A minimal Markov basis for $M'$ of the main effect model of a $2^{p-1}$
fractional factorial design of resolution $p$ is constructed as a set of 
square-free degree-$2$ moves.
\label{thm:main-theorem}
\end{theorem}

 In order to prove this theorem, we introduce the necessary notation.
We treat the frequency vector $\Bx \in \mathbb{N}^k$ as a 
$(p-1)$-dimensional contingency table, i.e., a $2^{p-1}$ table.
Denote the levels of each axis as $\{0,1\}$, and let
\[
{\cal I} = \{\Bi\} = \{(i_1\cdots i_{p-1})\} = \{0,1\}\times \cdots
\times \{0,1\}
\]
denote the set of cells. The frequency vector $\Bx \in \mathbb{N}^k$
is then written as
\[
\Bx = \{x(\Bi)\}_{\Bi \in {\cal I}} = \{x(i_1\ldots i_{p-1}),\ 
i_m \in \{0,1\},\ i = 1,\ldots,m\}.
\]
 Note that we express the $(p-1)$-dimensional cell $\Bi$ as a 
concatenation $i_1\cdots i_{p-1}$.
Since $\Bx$ has many zero entries, it is useful to 
express $\Bx$ with total $n$ as
$\Bx = \Bi_1\cdots\Bi_n$, where $\Bi_1,\ldots,\Bi_n$ are the
cells of positive frequencies of $\Bx$. For the case in which $x(\Bi) > 1$,
$\Bi$ is repeated $x(\Bi)$ times. Similarly, we express a move $\Bz$
as $\Bz = \Bi_1\cdots \Bi_n - \Bj_1\cdots \Bj_n$, where
$\Bi_1,\cdots, \Bi_n$ are the cells of positive elements of $\Bz$ and 
$\Bj_1,\cdots, \Bj_n$ are the cells of negative elements of $\Bz$.
We also write the positive and negative parts 
of $\Bz = \Bi_1\cdots \Bi_n - \Bj_1\cdots\Bj_n$ as
$\Bz^+ = \Bi_1\cdots\Bi_n$ and 
$\Bz^- = \Bj_1\cdots\Bj_n$, respectively.
 See Section 7.2 of \cite{Aoki-Hara-Takemura-2012} for the details of 
this notation.
Let ${\cal I} = {\cal I}_0 \cup {\cal I}_1$ and 
${\cal I}_0 \cap {\cal I}_1 = \emptyset$,  where
\[
{\cal I}_0 = \left\{(i_1\cdots i_{p-1})\ \left|\ 
\sum_{m=1}^{p-1}i_m = 0\ \ ({\rm mod}\ 2)\right.\right\},
\]
\[
{\cal I}_1 = \left\{(i_1\cdots i_{p-1})\ \left|\ 
\sum_{m=1}^{p-1}i_m = 1\ \ ({\rm mod}\ 2)\right.\right\}.
\]
We define the one-dimensional sum of $\Bx$, as
\[
x_m(i_m) = 
\sum_{i_1,\ldots,i_{m-1},i_{m+1},\ldots,i_{p-1}}
x(i_1\cdots i_{m-1}i_m i_{m+1}\cdots i_{p-1}).
\]
We also define the diagonal sum of $\Bx$ as
\[
x_D(0) = \sum_{\Bi \in {\cal I}_0}x(\Bi),\ \  
x_D(1) = \sum_{\Bi \in {\cal I}_1}x(\Bi).
\]
The sufficient statistic for the observation $\Bx$ is then written as
\begin{equation}
M'\Bx = (x_1(0),x_1(1),\ldots,x_{p-1}(0),x_{p-1}(1),x_D(0),x_D(1))'.
\label{eqn:sufficient-statistics-Mx}
\end{equation}
We also use this marginal notation for moves. Therefore, for a move
$\Bz$, we have $z_m(i_m) = 0$ for $m=1,\ldots,p-1,\ i_m = 0,1$ 
and $z_D(i_m) = 0$ for $i_m = 0,1$. We write levels $\{i_m, i_m^*\}$ to represent that $(i_m, i_m^*)$ is $(0,1)$ or $(1,0)$. 

\begin{example} 
Consider the case $p=4$. The set of cells 
 ${\cal I} = \{0,1\}\times \{0,1\}\times \{0,1\}$ is expressed as
${\cal I} = {\cal I}_0 \cup {\cal I}_1$, where
\[\begin{array}{c}
{\cal I}_0 = \{(000), (011), (110), (101)\},\\
{\cal I}_1 = \{(001), (010), (100), (111)\}.
\end{array}
\]
The sufficient statistic for the observation $\Bx = \{x(i_1i_2i_3)\}$ is
written as
\[
M'\Bx = (x_1(0),x_1(1),x_2(0),x_2(1),x_3(0),x_3(1),x_D(0),x_D(1))',
\]
where the one-dimensional sum of $\Bx$ is given as
\[
x_1(i_1) = \sum_{i_2}\sum_{i_3}x(i_1i_2i_3) = x(i_100)+x(i_101)+x(i_110)+x(i_111)\ \ \mbox{for}\ \ i_1 = 0,1,
\]
for example, and the diagonal sum of $\Bx$ is given as
\[\begin{array}{c}
x_D(0) = x(000)+x(011)+x(110)+x(101),\\ 
x_D(1) = x(001)+x(010)+x(100)+x(111). 
\end{array}
\]
\end{example}

We first  prove an important structure of moves for $M'$ for the 
case of $p=4$.

\begin{lemma}
Any $\Bz = \{z(i_1i_2i_3)\} \in \mathbb{Z}^{8}$ with 
either of the following types cannot be a move for $M'$ given by 
(\ref{eqn:M-p=4}).
\begin{enumerate}
\item $z(i_1i_2i_3),\ z(i_1^*i_2i_3),\ z(i_1^*i_2^*i_3),\ z(i_1^*i_2i_3^*) > 0$
for some $(i_1i_2i_3)$. 
\item $z(i_1i_2i_3),\ z(i_1i_2^*i_3),\ z(i_1^*i_2^*i_3),\ z(i_1i_2^*i_3^*) > 0$
for some $(i_1i_2i_3)$. 
\item $z(i_1i_2i_3),\ z(i_1i_2i_3^*),\ z(i_1^*i_2i_3^*),\ z(i_1i_2^*i_3^*) > 0$
for some $(i_1i_2i_3)$. 
\end{enumerate}
\label{lemma-p4-lemma}
\end{lemma}
\paragraph*{Proof.}\  
Consider the case 1. In this case, we have
\[
z(i_1i_2i_3^*),\ z(i_1i_2^*i_3),\ z(i_1i_2^*i_3^*),\ z(i_1^*i_2^*i_3^*) < 0 
\]
From $M'\Bz = \Bzero$. Because the one-dimensional sum $z_1(i_1)$ is zero, 
we have
\[
z(i_1i_2i_3) = |z(i_1i_2i_3^*)|+ |z(i_1i_2^*i_3)|+ |z(i_1i_2^*i_3^*)|.
\]
On the other hand, because the one-dimensional sum $z_2(i_1)$ is zero,
we also have
\[
z(i_1i_2i_3) + z(i_1^*i_2i_3)+ z(i_1^*i_2i_3^*) = |z(i_1i_2^*i_3^*)|.
\]
From these relations, we have
\[
|z(i_1i_2i_3^*)|+ |z(i_1i_2^*i_3)|+ z(i_1^*i_2i_3)+ z(i_1^*i_2i_3^*) = 0,
\]
which contradicts the assumption. The cases 2 and 3 are also proved
 similarly.
\hspace*{\fill}$\Box$

\vspace*{10mm}

Next, we  prove Theorem \ref{thm:main-theorem}. 
In the following proof, we use 
the method of distance reduction.  
The theorem is proved
by contradiction. First, we assume that the target set of moves, i.e., 
the set of square-free degree-$2$ moves, is not a Markov basis, and
derive a contradiction by showing that the distance of any two states 
in the same fiber can be reduced by the moves from the target set.
See Chapter 6 of \cite{Aoki-Hara-Takemura-2012}, for details.

\paragraph*{Proof of Theorem \ref{thm:main-theorem}.}\ 
Write the set of square-free degree-$2$ moves for $M'$ as ${\cal B}^*$. 
In order to  prove that ${\cal B}^*$ is a Markov basis for $M'$, 
assume that ${\cal B}^*$ is not a Markov basis. 
Then, for some $\By^o \in \mathbb{N}^k$, 
${\cal F}(M'\By^o)$ consists of more than one ${\cal B}^*$-equivalence 
class. 
Let ${\cal H}_1$ and ${\cal H}_2$ denote two different ${\cal B}^*$-equivalence
classes. Choose $\Bx \in {\cal H}_1$ and $\By\in {\cal H}_2$ such that
\[
|\Bz| = |\By - \Bx| = \sum_{\Bi \in {\cal I}}|y(\Bi) - x(\Bi)|
\]
is minimized. Since $\Bx$ and $\By$ are chosen from different ${\cal B}^*$ 
classes, this minimum must be positive. 
We reveal a contradiction by induction on $p$.

Consider the case of $p=4$.  Assume without loss of generality 
that $z(000) > 0$.
Here, we consider $z(111)$. 
If $z(111) > 0$, then $y(000), y(111) > 0$ follows. 
Since $\Bz$ is a move, $\Bz$ also 
has negative entries, $x(i_1i_2i_3) < 0$ for some $(i_1i_2i_3)$. Then, we can
subtract the move
\[
\Bz_1 = (000)(111) - (i_1i_2i_3)(i_1^*i_2^*i_3^*) \in {\cal B}^*
\]
from $\By$ without generating negative entries. However, we have
$|(\By - \Bz_1) - \Bx| < |\By - \Bx|$ 
because $\Bz_1$ and $\Bx$ have a common support, which contradicts the 
definition of $\Bx, \By$. 
Then, we have $z(111) \leq 0$. 
Next, we assume that $z(001) < 0$. In this case, by symmetry (i.e., 
for same reason described above),
we can derive a contradiction if $z(110)<0$. Then, we have $z(110) \geq 0$.
Next, we consider $z(010)$. 
\begin{enumerate}
\item Assume $z(010) < 0$. 
By symmetry, we have $z(101) \geq 0$. 
From $z_D(0)=0$, we have $z(011) < 0$. 
From $z_D(1)=0$, we have $z(100) > 0$. 
From $z_1(1)=0$, we have $z(111) < 0$. 
However, $z(001), z(010), z(011), z(111) < 0$ contradicts the idea that $-\Bz$
is a move for $M'$ as stipulated by Lemma \ref{lemma-p4-lemma}.
\item Assume that $z(010) \geq 0$.
From $z_3(0)=0$, we have $z(100)<0$. 
By symmetry, we have $z(011) \geq 0$. 
From $z_2(1)=0$, we have $z(111)<0$. 
From $z_D(0)=0$, we have $z(101)<0$.
However, $z(001), z(100), z(101), z(111) < 0$ contradicts the idea that $-\Bz$
is a move for $M'$ as stipulated by Lemma \ref{lemma-p4-lemma}.
\end{enumerate}
Then, we have  proved  the theorem for $p=4$.

Next, we assume that $p\geq 5$ and that the theorem holds for cases up to $p-1$. 
We begin with the case in which some values of $M'\By^o$ are zero. 
Without loss of generality, we have two cases.
\begin{enumerate}
\item In the case of $y_{p-1}^o(1) = 0$, define the set
\begin{equation}
\left\{\By \ \left|\ 
y(i_1\ldots i_{p-2}) = x(i_1\cdots i_{p-2}0) 
\ {\rm for}\ \Bx \in {\cal F}(M'\By^o)
\right.\right\}.
\label{eq:0-case-1}
\end{equation}
Then, we have 
\[
\left\{
\begin{array}{l}
y_m(i_m) = x_m(i_m),\ m=1,\ldots,p-2,\ i_m = 0,1,\\
y_D(i_m) = x_D(i_m),\ i_m = 0,1.
\end{array}
\right.
\]
Therefore, from the definition of the fiber, the set (\ref{eq:0-case-1}) 
coincides with the fiber for
\[
(x_1(0),x_1(1),\ldots,x_{p-2}(0),x_{p-2}(1),x_D(0),x_D(1))'.
\]
Based on the assumption of induction, this fiber is connected by the 
set of square-free degree-$2$ moves ${\cal B} \subset \mathbb{Z}^{2^{p-2}}$.
Then, the set of square-free degree-$2$ moves 
\[
\left\{
\widetilde{\Bz} = \{\widetilde{z}(\Bi)\}\ \left|\ 
\widetilde{z}(i_1\cdots i_{p-2}i_{p-1}) = \left\{
\begin{array}{ll}
z(i_1\cdots i_{p-2}), & {\rm if}\ i_{p-1}=0\\
0, & {\rm if}\ i_{p-1}=1
\end{array}
\right.\ \ {\rm for}\ \Bz \in {\cal B}
\right.
\right\}
\]
connects the fiber for $M'\By^o$.
\item In the case of $y_D^o(1) = 0$, define the set 

\begin{equation}
\left\{\By \ \left|\ 
\begin{array}{l}
{
y(i_1\ldots i_{p-3}0) = \left\{\begin{array}{ll}
x(i_1\cdots i_{p-3}00), & {\rm if}\ i_1+\cdots i_{p-3}=0\ ({\rm mod}\ 2)\\ 
x(i_1\cdots i_{p-3}01), & {\rm if}\ i_1+\cdots i_{p-3}=1\ ({\rm mod}\ 2)
\end{array}\right.
}\\
{
y(i_1\ldots i_{p-3}1) = \left\{\begin{array}{ll}
x(i_1\cdots i_{p-3}11), & {\rm if}\ i_1+\cdots i_{p-3}=0\ ({\rm mod}\ 2)\\ 
x(i_1\cdots i_{p-3}10), & {\rm if}\ i_1+\cdots i_{p-3}=1\ ({\rm mod}\ 2)
\end{array}\right.
}
\end{array}
{\rm for}\ \Bx \in {\cal F}(M'\By^o)
\right.\right\}.
\label{eq:0-case-2}
\end{equation}
Then, we have
\[
\left\{
\begin{array}{l}
y_m(i_m) = x_m(i_m),\ m=1,\ldots,p-2,\ i_m = 0,1\\
y_{p-2}(i_m) = x_{p-2}(i_m),\ i_m = 0,1\\
y_D(i_m) = x_{p-1}(i_m),\ i_m = 0,1. 
\end{array}
\right.
\]
Therefore, based on the assumption of induction, we can similarly  prove that 
the set of square-free degree-$2$ moves connects the fiber for $M'\By^o$. 
\end{enumerate}

Based on the above consideration, we can restrict our attention to the cases
in which all of the elements of $M'\By^o$ are positive. 
Assume that ${\cal F}(M'\By^o)$ is such a fiber, and return to the consideration 
of $\Bz = \By - \Bx$ for $\Bx \in {\cal H}_1$, $\By \in {\cal H}_2$, where
${\cal H}_1$ and ${\cal H}_2$ are two different ${\cal B}^*$-equivalence
classes of ${\cal F}(M'\By^o)$. 
We must consider the cases in which $p\ (\geq 5)$ is either even or odd.

Assume that $p$ is even. 
Since we assume that $|\Bz|$ is positive, we can assume without loss of generality that $z(0\cdots 0) > 0$.
Since all of the elements of $M'\By^o$ are positive, there is at least 
one more positive entry in $\Bz$. If $z(1\cdots 1) > 0$, we can obtain 
$z(i_1\cdots i_{p-1}) < 0$ for some $(i_1\cdots i_{p-1}) \in {\cal I}_0$.
Then, we can subtract a move
\[
\Bz_1 = (0\cdots 0)(1\cdots 1) - (i_1\cdots i_{p-1})(i_1^*\cdots i_{p-1}^*)
\in {\cal B}^* 
\]
from $\By$ without generating negative entries, 
and we have $|(\By - \Bz_1) - \Bx| < |\By - \Bx|$,
which is a contradiction. Then, we can assume $z(0\cdots 0) \leq 0$, and
we can choose other positive entries in $\Bz$.
Without loss of generality, we can assume
\[
z(\underbrace{0\cdots 0}_{p-1-r}\underbrace{1\cdots 1}_{r}) > 0,
\]
where $r$ is odd and $p-1-r \geq 2$ is even. 
Since the one-dimensional totals $y_1^o(0),\ldots,y_{p-1-r}^o(0)$ are 
positive, $\Bz$ has another positive element. 
Assume that $\Bz$ has a positive element expressed as
\[
z(\underbrace{\ast\cdots\ast 1\ast\cdots\ast 1\ast\cdots\ast}_{p-1-r}\underbrace{\ast\cdots \ast}_{r}) > 0,
\]
 where we express $\ast$ as $0$ or $1$. 
Rewrite this as 
\[
z(\underbrace{\ast\cdots \ast}_{p-3-r}11\underbrace{\ast\cdots \ast}_{r}) > 0
\]
without loss of generality. Then, we can subtract a move
\[
\Bz_2 = 
(\underbrace{0\cdots 0}_{p-3-r}00\underbrace{1\cdots 1}_{r})
(\underbrace{\ast\cdots \ast}_{p-3-r}11\underbrace{\ast\cdots \ast}_{r})
-
(\underbrace{0\cdots 0}_{p-3-r}11\underbrace{1\cdots 1}_{r})
(\underbrace{\ast\cdots \ast}_{p-3-r}00\underbrace{\ast\cdots \ast}_{r})
\in {\cal B}^*
\]
from $\By$ without generating negative entries. 
Note that $|(\By - \Bz_2) - \Bx|$ is not necessarily smaller than 
$|\By - \Bx|$. However, 
$(\By - \Bz_2) - \Bx = \Bz - \Bz_2 (= \Bz^*)$ has positive entries
\[
z^*(0\cdots 0) > 0,\ 
z^*(\underbrace{0\cdots 0}_{p-3-r}\underbrace{1\cdots 1}_{r}) > 0.
\]
Therefore, repeating this procedure, we can reduce 
the pattern of the positive part of $\Bz$ as
\[
(0\cdots 0)
(\underbrace{0\cdots 0}_{p-1-r}\underbrace{1\cdots 1}_{r})
(\underbrace{10\cdots 0}_{p-1-r}\underbrace{\ast\cdots \ast}_{r})
(\underbrace{010\cdots 0}_{p-1-r}\underbrace{\ast\cdots \ast}_{r})
\cdots
(\underbrace{0\cdots 01}_{p-1-r}\underbrace{\ast\cdots \ast}_{r}).
\]
Now, consider the cells
\[
(\underbrace{10\cdots 0}_{p-1-r}\underbrace{\ast\cdots \ast}_{r}),
(\underbrace{010\cdots 0}_{p-1-r}\underbrace{\ast\cdots \ast}_{r}),
\cdots, 
(\underbrace{0\cdots 01}_{p-1-r}\underbrace{\ast\cdots \ast}_{r}).
\]
If $p-1-r > 2$, there are two cells
in these $p-1-r$ cells, both of which have the same level 
in the same axis in the latter $r$ axes, e.g.,
\[
(\underbrace{10\cdots 0}_{p-1-r}\underbrace{0\ast\cdots \ast}_{r}),
(\underbrace{010\cdots 0}_{p-1-r}\underbrace{0\ast\cdots \ast}_{r})
\]
or
\[
(\underbrace{10\cdots 0}_{p-1-r}\underbrace{1\ast\cdots \ast}_{r}),
(\underbrace{010\cdots 0}_{p-1-r}\underbrace{1\ast\cdots \ast}_{r}).
\]
In the former case, 
the three-dimensional marginal table of $\Bz^+$ with respect to 
$i_1, i_2, i_{p-r}$ has positive elements
\[
(000)(001)(100)(010),
\]
which contradicts Lemma \ref{lemma-p4-lemma}.
Similarly, in the latter case, 
the three-dimensional marginal table of $\Bz^+$ with respect to 
$i_1, i_2, i_{p-r}$ has positive elements
\[
(000)(001)(101)(011),
\]
which also contradicts Lemma \ref{lemma-p4-lemma}.
Therefore, we need only consider the case in which $p-1-r = 2$ and
the positive entries of $\Bz$ are written as
\[
(000\cdots 0)
(001\cdots 1)
(10i_3 \cdots i_{p-1})
(01i_3^* \cdots i_{p-1}^*).
\]
In this case, we can subtract a move
\[
(10i_3 \cdots i_{p-1})
(01i_3^* \cdots i_{p-1}^*)
-
(100\cdots 0)
(011\cdots 1) \in {\cal B}^*
\]
from $\By$ without generating negative entries, and we have
positive elements
\[
(000\cdots 0)
(001\cdots 1)
(100\cdots 0)
(011\cdots 1).
\]
Here, considering the collapsing table with respect to $i_3,\ldots,i_{p-1}$ axes,
such as $(0\cdots 0) \rightarrow 0$ and 
$(1\cdots 1) \rightarrow 1$, we again have the $p=4$ case and can derive
a contradiction from Lemma \ref{lemma-p4-lemma}. 
Therefore, we have proved the theorem for the case in which $p$ is even.

On the other hand, for the case in which $p$ is odd, 
we focus on $(01\cdots 1)$ 
rather than $(1\cdots 1)$ and can derive the proof in a similar manner.
Thus, the proof of the theorem is complete.
\hspace*{\fill}$\Box$

\subsection{Structure of the minimal Markov basis}
Based on the proof of Theorem \ref{thm:main-theorem}, the structure of
the minimal Markov basis is revealed 
for the main effect models of
$2^{p-1}$ fractional factorial designs of resolution $p$. 
We refer to a fiber ${\cal F}(M'\By^o)$ as a degree-$2$ fiber if 
$\sum_{\Bi \in {\cal I}} y^o(\Bi) = 2$. 
As a consequence of Theorem \ref{thm:main-theorem}, we see
that the degree-$2$ fiber plays an important role 
in the construction of a Markov basis.

\begin{corollary}
A minimal Markov basis 
for $M'$ of the main effect models of
$2^{p-1}$ fractional factorial designs of resolution $p$ 
is constructed as the set of
moves connecting all degree-$2$ fibers with more than one element.
\end{corollary}

Based on this corollary, we can characterize the structure of the minimal 
Markov basis considering all degree-$2$ fibers with more than
one element. We refer to such a fiber as an {\it essential fiber}.
 We consider the essential fibers of small $p$ below. 
Recall that the sufficient statistic is given as
(\ref{eqn:sufficient-statistics-Mx}).
In the following, deleting commas from 
(\ref{eqn:sufficient-statistics-Mx}), 
we write the sufficient statistic as
\[
M'\Bx = (x_1(0)x_1(1),\ldots,x_{p-1}(0)x_{p-1}(1),x_D(0)x_D(1))'.
\]

\begin{enumerate}
\item  Case $p=4$. There is only one essential fiber, which is
given by
\[
{\cal F}(11,11,11,11) = \{(111)(222),\ (112)(221),\ (121)(212),\ 
(122)(211)\}.
\]
In order to connect this $4$-element fiber, $3$ moves,  for example 
\[
\{
(111)(222) - (112)(221),\ 
(111)(222) - (121)(212),\ 
(111)(222) - (122)(211)\},
\]
constitute a minimal Markov basis.

\item  Case $p=5$. There are $10$ essential $4$-element 
fibers, e.g.,
\[
{\cal F}(20,11,11,11,11),\  
{\cal F}(02,11,11,11,11),\  
\ldots,
{\cal F}(11,11,11,11,20),\  
{\cal F}(11,11,11,11,02).
\]
Connecting these $4$-element fibers requires $3$ moves for each fiber.
Therefore, the minimal Markov basis is constructed as $3\times 10 = 30$ moves.

\item  Case $p=6$. There are $60$ essential $4$-element 
fibers, e.g.,
\[
{\cal F}(20,20,11,11,11,11),\  
\ldots,
{\cal F}(11,11,11,11,02,02)
\]
and one essential $16$-element fiber
\[
{\cal F}(11,11,11,11,11,11).
\]
Connecting these fibers requires $3$ moves for each $4$-element
fiber and $15$ moves for the $16$-element fiber. 
Therefore, the minimal Markov basis is constructed 
of $3\times 60 + 15 = 195$ moves.

\item  Case $p=7$. There are $280$ essential $4$-element 
fibers, e.g.,
\[
{\cal F}(20,20,20,11,11,11,11),\  
\ldots,
{\cal F}(11,11,11,11,02,02,02)
\]
and $14$ essential $16$-element fibers,
e.g.,
\[
{\cal F}(20,11,11,11,11,11,11),\  
\ldots,
{\cal F}(11,11,11,11,11,11,02).
\]
The minimal Markov basis is constructed of $3\times 280 + 15\times 14 = 1,050$ 
moves.

\item  Case $p=8$. There are $1,120$ essential $4$-element 
fibers, e.g.,
\[
{\cal F}(20,20,20,20,11,11,11,11),\  
\ldots,
{\cal F}(11,11,11,11,02,02,02,02),
\]
$112$ essential $16$-element fiber, e.g.,
\[
{\cal F}(20,20,11,11,11,11,11,11),\  
\ldots,
{\cal F}(11,11,11,11,11,11,02,02)
\]
and
one essential $64$-element fiber
\[
{\cal F}(11,11,11,11,11,11,11,11).
\]
The minimal Markov basis is constructed as 
$3\times 1120 + 15\times 112 + 63 = 5,103$ moves.
\end{enumerate} 

We summarize the number of essential fibers as follows.
\begin{proposition}
A minimal Markov basis 
for $M'$ of the main effect models of
$2^{p-1}$ fractional factorial design 
of resolution $p$ 
is constructed as 
\[
\displaystyle\bigcup_{t = 1}^{s - 1}
\left\{{\rm moves\ connecting}\ 
4^{t-1}{p \displaystyle\choose{2(t-1)}}\ {\rm essential}\ 
\mbox{$4^{s-t}$-{\rm element\ fibers}} 
\right\}
\]
for $p=2s$, and 
\[
\displaystyle\bigcup_{t = 1}^{s - 1}
\left\{{\rm moves\ connecting}\ 
2^{2t-1}{p \displaystyle\choose{2t-1}}\ {\rm essential}\ 
\mbox{$4^{s-t}$-{\rm element\ fibers}} 
\right\}
\]
for $p = 2s+1$.
\end{proposition}

\subsection{Example} 
\begin{figure*}
\begin{center}
\includegraphics[width=90mm, height=60mm]{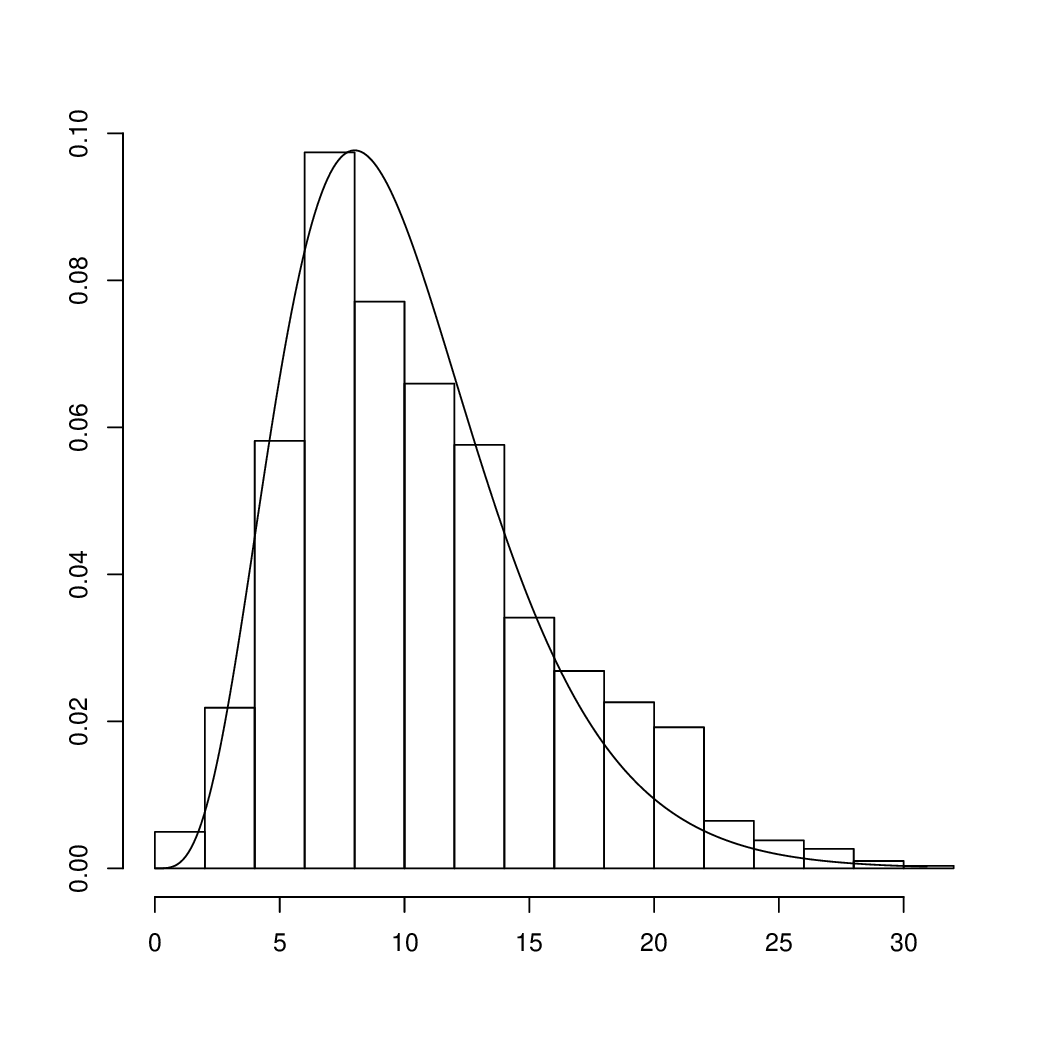}
\end{center}
\caption{Asymptotic and Monte Carlo estimated distribution of the
 likelihood ratio statistic.}
\label{fig:hist-asymp}
\end{figure*}

To illustrate our method, we show an example of analyzing actual data.
Table \ref{tbl:chemical-reaction} shows data of chemical reaction from 
Section 9 of \cite{onyiah-book}.
The five two-level 
factors, which are thought to influence the yield of a chemical
reaction in this experiment, are 
$\BA$:\ Temperature ($30^{\circ}$C or $40^{\circ}$C), 
$\BB$:\ Pressure ($1$ atm or $2$ atm), 
$\BC$:\ Reaction time ($2$ hours or $4$ hours), 
$\BD$:\ Reactant concentration ($30$ \% or $40$ \%) and
$\BE$:\ Catalyst (Absent or Present). 
For each level of factors in the $2^{5-1}_{{\rm V}}$ 
fractional factorial design, the value of the resulting responses ($\By$)
is reported.
\begin{table*}
\label{tbl:chemical-reaction}
\begin{center}
\caption{Design and resulting responses in the chemical reaction data}
\begin{tabular}{rrrrrrcc}\hline
& \multicolumn{5}{c}{Factor} & $\By$ & $\By$ (rounded)\\
Run & $\BA$ & $\BB$ & $\BC$ & $\BD$ & $\BE$ & \\ \hline
$1$ & $1$ &  $1$ &  $1$ &  $1$ &  $1$ & $30.3$ &  $30$\\ 
$2$ &$-1$ &  $1$ &  $1$ &  $1$ & $-1$ & $27.6$ & $28$\\
$3$ & $1$ & $-1$ &  $1$ &  $1$ & $-1$ & $25.8$ &  $26$\\
$4$ & $-1$ & $-1$ &  $1$ &  $1$ &  $1$ & $31.3$ & $31$\\
$5$ & $1$ &  $1$ & $-1$ &  $1$ & $-1$ &  $24.9$ & $25$\\
$6$ & $-1$ &  $1$ & $-1$ &  $1$ &  $1$ & $36.3$ & $36$\\
$7$ &  $1$ & $-1$ & $-1$ &  $1$ &  $1$ & $24.9$ & $25$\\
$8$ & $-1$ & $-1$ & $-1$ &  $1$ & $-1$ & $28.3$ & $28$\\
$9$ & $1$ &  $1$ &  $1$ & $-1$ & $-1$ &  $37.9$ & $38$\\
$10$ &$-1$ &  $1$ &  $1$ & $-1$ &  $1$ &  $34.9$ & $35$\\
$11$ &  $1$ & $-1$ &  $1$ & $-1$ &  $1$ & $35.9$ & $36$\\
$12$ &$-1$ & $-1$ &  $1$ & $-1$ & $-1$ &  $28.3$ & $28$\\
$13$ &  $1$ &  $1$ & $-1$ & $-1$ &  $1$ & $33.6$ & $34$\\
$14$ & $-1$ &  $1$ & $-1$ & $-1$ & $-1$ & $37.0$ & $37$\\
$15$ &  $1$ & $-1$ & $-1$ & $-1$ & $-1$ & $32.2$ & $32$\\
$16$ & $-1$ & $-1$ & $-1$ & $-1$ &  $1$ & $37.4$ & $37$\\ \hline
\end{tabular}
\end{center}
\end{table*}
Though the original response $\By$ is nonnegative real value vector, 
we round $\By$
to nonnegative integer and treat $\By$ as a realisation of Poisson 
random vector as described in Section 2. Though it 
loses some information of data, we round $\By$ to use our method. 
The test statistic we use is the likelihood ratio statistic
(\ref{eq:g2}), where the fitted value under the main effect model,
$\{\mu_i\}$, is 
calculated as
\[\begin{array}{cl}
( & 35.80,\ 27.13,\ 31.05,\ 32.04,\ 33.61,\ 35.52,\ 29.37,\ 28.01,\\
  & 36.67,\ 29.13,\ 25.47,\ 27.78,\ 33.34,\ 38.45,\ 32.56,\ 30.07\ \ ).
\end{array}
\]
The value of the likelihood ratio statistic is $3.37$ with $10$ degree
of freedom. The asymptotic $p$ value based on $\chi^2_{10}$ is $0.9712$. 
To estimate $p$ value by the Markov chain Monte Carlo method, first we
need to calculate a Markov basis. From the result of Section 3.3, the
minimal Markov basis is constructed by $30$ square-free degree $2$
moves. Using this minimal Markov basis, 
$p$ value is estimated as $0.96$ from 
$n = 10000$ samples after $10000$ burn-in samples. 
These $p$ values suggest that the fitting of the main effect model to
the data is fine.
Figure \ref{fig:hist-asymp} shows a histogram of the Monte Carlo sampling
generated from the conditional distribution of the likelihood ratio
statistic under the main effect model, along with the corresponding
asymptotic distribution $\chi^2_{10}$.

\section{Discussion}
In the present paper, we consider Markov chain Monte Carlo tests for the
factor effects in designed experiments. 
As noted in Section 1, 
the Markov chain Monte Carlo procedure is a valuable tool when the exact 
calculation
of the $p$ value is infeasible. 
In addition, the closed-form expression of the normalizing constant for 
the null distribution of our model is not obtained.
Therefore, the Markov chain Monte Carlo procedure is  a valuable 
alternative.

In order to perform Markov chain Monte Carlo tests, it is often difficult 
to calculate a Markov 
basis. As reported in previous studies, the structure of the Markov basis
is generally complicated. On the other hand, problems in the Markov 
basis that have a simple structure are also important. In particular, 
special cases
in which a Markov basis, i.e., the generator of the toric ideal, 
is constructed only by degree-$2$ moves have attracted the attention of 
algebraists.  See \cite{Ohsugi-Hibi-1999}, for example.
The present paper demonstrates that the problem considered 
 here is one such important case. 

 Though our result in this paper is restricted to the
main effect model of 
$1/2$-fractional factorial designs, it is attractive topic 
to use the idea of this paper for
$1/4$- or $1/8$-fractional factorial designs.
 It is more important to consider models that include various 
interaction effects.
A general method by which to treat interaction effects is presented 
in \cite{Aoki-Takemura-2010}. However, the structure of the Markov basis
for models including interaction effects generally becomes complicated.
In addition, the Markov chain Monte Carlo approach can also be used for
general fractional factorial designs. Because the structure of 
the Markov basis in general settings is intractable, it is important to
also consider well-known designs such as Plackett-Burman designs or balanced 
incomplete block designs.

\bibliographystyle{plain}

\end{document}